\documentclass[12pt]{amsart}
\usepackage{}

\usepackage{amsmath}
\usepackage{amsfonts}
\usepackage{amssymb}
\usepackage[all]{xy}           

\usepackage{bbding}
\usepackage{txfonts}
\usepackage{amscd}

\usepackage[shortlabels]{enumitem}
\usepackage{ifpdf}
\ifpdf
  \usepackage[colorlinks,final,backref=page,hyperindex]{hyperref}
\else
  \usepackage[colorlinks,final,backref=page,hyperindex]{hyperref}
\fi
\usepackage{tikz}
\usepackage[active]{srcltx}

\makeatletter

\newtheorem{thm}{Theorem}[section]
\newtheorem{lem}[thm]{Lemma}
\newtheorem{cor}[thm]{Corollary}
\newtheorem{pro}[thm]{Proposition}
\newtheorem{ex}[thm]{Example}

\newtheorem{defi}[thm]{Definition}

\newcommand {\emptycomment}[1]{}

\newcommand{\cb}{\mathbb C}
\newcommand{\HH}{\mathbb H}

\newcommand{\Der}{\mathrm{Der}}

\newcommand{\Inn}{\mathrm{Inn}}

\newcommand{\Img}{\mathrm{Im}}

\def\oll#1#2{\mathrel {{\triangleright}_{(#1)}^{#2}}}
\def\orr#1#2{\mathrel {{\triangleleft}_{(#1)}^{#2}}}

\def\ool#1{\mathrel {{\triangleright}_{(#1)}}}

\def\oor#1{\mathrel {{\triangleleft}_{(#1)}}}
\def\oo#1{\mathrel {{}_{(#1)}}}
\def\ooo#1#2{\mathrel {{}_{(#1)}^{#2}}}

\def\id{\mathop {\fam0 id}\nolimits}

\def\Ker{\mathrm{Ker}\,}
\def\Hom{\mathrm{Hom}}

\def\HH{\mathrm{HH}}

\begin{document}

\title[Gerstenhaber algebra of an associative conformal algebra]
{Gerstenhaber algebra of the Hochschild cohomology of an associative
conformal algebra}

\author{Bo Hou}
\address{School of Mathematics and Statistics, Henan University, Kaifeng 475004,
China}
\email{bohou1981@163.com}
\vspace{-5mm}

\author{Zhongxi Shen}
\address{School of Mathematics and Statistics, Henan University, Kaifeng 475004,
China}
\email{18738951378@163.com}
\vspace{-5mm}

\author{Jun Zhao}
\address{School of Mathematics and Statistics, Henan University, Kaifeng 475004,
China}
\email{zhaoj@henu.edu.cn}


\begin{abstract}
We define a cup product on the Hochschild cohomology of an associative conformal
algebra $A$, and show the cup product is graded commutative. We define a graded Lie
bracket with the degree $-1$ on the Hochschild cohomology $\HH^{\ast}(A)$ of an
associative conformal algebra $A$, and show that the Lie bracket together with
the cup product is a Gerstenhaber algebra on the Hochschild cohomology of an
associative conformal algebra.
Moreover, we consider the Hochschild cohomology of split extension
conformal algebra $A\widehat{\oplus}M$ of $A$ with a conformal bimodule $M$, and
show that there exists an algebra homomorphism from $\HH^{\ast}(A\widehat{\oplus}M)$
to $\HH^{\ast}(A)$.
 \end{abstract}


\maketitle




\vspace{-4mm}

\section{Introduction}\label{sec:intr}

The notion of a conformal algebra was introduced by Kac in 1997, which encodes an axiomatic
description of the operator product expansion of chiral fields in conformal field
theory. The theory of Lie conformal algebras appeared as a formal language
describing the algebraic properties of the operator product
expansion in two-dimensional conformal field theory (\cite{BPZ,Ka,Ka1}) and it has
received a lot of attention.
Structure theory and representation theory of Lie conformal algebras are
widely studied in a number of publications (see \cite{BKV,CK,CKW,DK}, and so on).

Naturally, associative conformal algebras come from representations
of Lie conformal algebras. Moreover, some Lie conformal algebras appeared
in physics are embeddable into associative one \cite{Ro, Ro1}.
The structure theory and representation theory of associative conformal algebras have
attracted the attention of many scholars and achieved a series of results, see
\cite{BFK,BFK1,BFK2,Ko,Ko1,Ko2,Ko3,Re,Re1}.
In the recent years there is a great advance in the cohomology theory
of associative conformal algebras. In \cite{Do} and \cite{Do1},
Dolguntseva introduced the Hochschild cohomology of associative conformal
algebras, and proved that the second Hochschild cohomology group of the
conformal Weyl algebra, ${\rm Cend}_{n}$ and ${\rm Cur}_{n}$ with values in
any bimodule is trivial. Kozlov showed that the second Hochshild cohomology group
of the associative conformal algebra ${\rm Cend}_{1,x}$ with values in any
bimodule is also trivial \cite{Koz}. Recently, Kolesnikov and Kozlov proved
that all semisimple algebras of conformal endomorphisms which have
the trivial second Hochschild cohomology group with coefficients in every
conformal bimodule (see \cite{KK}). In \cite{Lib}, Liberati defined the
cohomology of associative H-pseudoalgebras and described the zero, the first
and the second cohomologies of associative H-pseudoalgebras.

In this paper, we consider the Gerstenhaber algebra structure on the
Hochschild cohomology of an associative conformal algebra.
It is a classical fact that the Hochschild cohomology $\HH^\bullet(A, A)$ of
an associative algebra $A$ carries a Gerstenhaber algebra structure \cite{Ger}.
A Gerstenhaber algebra is a graded commutative, associative algebra
$\Big(\bigoplus_{i\in\mathbb{Z}} \mathcal{A}^i,\, \sqcup\Big)$ together with
a degree $-1$ graded Lie bracket $[-,-]$ compatible with the product $\sqcup$
in the sense of the following Leibniz rule
\begin{align}\label{g-com}
[a\sqcup b, c] = [a, c]\sqcup b + (-1)^{(|c|-1) |a|} a\sqcup [b, c].
\end{align}
For an associative algebra $A$, the Hochschild cochain complex
$\mathcal{C}^{\ast}(A, A)$
carries a cup product $\sqcup: \mathcal{C}^{m}(A, A)\times \mathcal{C}^{n}(A, A)
\rightarrow \mathcal{C}^{m+n}(A, A)$
defined by
\begin{align*}
(f\sqcup g)(a_{1},\dots, a_{m+n})
=f(a_{1}, \dots, a_{m})g(a_{m+1},\dots, a_{m+n}).
\end{align*}
It turns out that the Hochschild coboundary operator $\delta$ is a graded derivation
with respect to the cup product. Hence, it induces a cup product $\sqcup$ on the
Hochschild cohomology $\HH^{\ast}(A, A)$. Moreover, the cochain groups
$\mathcal{C}^{\ast}(A, A)$ carry a degree $-1$ graded Lie bracket compatible with
the Hochschild cobundary \cite{Ger}. Therefore, it gives rise to a
degree $-1$ graded Lie bracket on $\HH^{\ast}(A, A)=\bigoplus_{i\geq0}\HH^{i}(A, A)$.
The cup product and the degree $-1$ graded Lie bracket on the Hochschild cohomology
$\HH^{\ast}(A, A)$ are compatible in the sense of (\ref{g-com}) to make it into a
Gerstenhaber algebra. Recently, for the Hom-associative, the Gerstenhaber algebra
structure on the Hochschild cohomology is given \cite{Das}.
Here we imitate the case of associative algebra, and consider
the Gerstenhaber algebra structure on the Hochschild cohomology of an associative
conformal algebra. In 1999, Bakalov, Kac and Voronov have given the definitions of
Hochschild cohomology complex and the Hochschild cohomology groups of associative
conformal algebras \cite{BKV}. Here we define a cup product and a Gerstenhaber
bracket on this complex, and show that these induces a Gerstenhaber algebra
structure on the Hochschild cohomology of an associative conformal algebra.

The paper is organized as follows.
In Section \ref{sec:prel}, we recall the notions of associative conformal
algebras, conformal bimodules, and the Hochschild cohomology of associative conformal
algebras. In Section \ref{sec:lie}, we define a Gerstenhaber bracket on the
Hochschild cohomology complex of an associative conformal algebra $A$. This bracket
induces a graded Lie bracket $[-,-]$ on $\HH^{\ast}(A, A)=\bigoplus_{i\geq0}
\HH^{i}(A, A)$ with degree $-1$. In Section \ref{sec:cup},
we define a product on the Hochschild cohomology complex of an associative
conformal algebra $A$. This product induces a cup product $\sqcup$
on the Hochschild cohomology $\HH^{\ast}(A, A)=\bigoplus_{i\geq0}\HH^{i}(A, A)$.
In Section \ref{sec:Gerst}, we show that $\HH^{\ast}(A, A)=
\bigoplus_{i\geq0}\HH^{i}(A, A)$ with the cup product $\sqcup$ and the graded Lie
bracket $[-,-]$ of degree $-1$ is a Gerstenhaber algebra.
In Section \ref{sec:exten}, we define the split extension associative conformal
algebra $A\widehat{\oplus}M$ of an associative conformal algebra $A$ by a conformal
bimodule $M$, consider the relationship between $\HH^{\ast}(A, A)$ and
$\HH^{\ast}(A\widehat{\oplus}M, A\widehat{\oplus}M)$, and give an algebra
homomorphism from $\HH^{\ast}(A\widehat{\oplus}M)$ to $\HH^{\ast}(A)$.

Throughout this note, we fix $\cb$ an algebraical closed field and characteristic zero
(for example, the field of complex numbers), denote by $\mathbb{Z}_{+}$ the set of all
nonnegative integers. All vector spaces are $\cb$-vector spaces, all linear maps and
bilinear maps are $\cb$-linear, all tensor products are over $\cb$, unless otherwise
specified. For any vector space $V$ and variable $\lambda$, we use $V[\lambda]$ to
denote the set of polynomials of $\lambda$ with coefficients in $V$.

\section{Associative conformal algebras and Hochschild cohomology} \label{sec:prel}

We recall the notions of associative conformal algebras, conformal bimodules over
associative conformal algebras, and the Hochschild cohomology of an associative
conformal algebra with coefficients in a bimodule. For the details see
\cite{BDK,DK,BKV}.

\begin{defi}\label{de:alg}
A {\rm conformal algebra} $A$ is a $\cb[\partial]$-module endowed with a
$\cb$-bilinear map $\cdot\oo{\lambda}\cdot : A\times A\rightarrow A[\lambda]$,
$(a, b)\mapsto a\oo{\lambda} b$ satisfying
\begin{eqnarray*}
\partial a\oo{\lambda}b=-\lambda a\oo{\lambda}b,  \qquad
a\oo{\lambda}\partial b=(\partial+\lambda)a\oo{\lambda}b,
\quad\forall\, a,b\in A.
\end{eqnarray*}
An {\rm associative conformal algebra} $A$ is a conformal algebra satisfying
\begin{eqnarray*}
(a\oo{\lambda}b)\oo{\lambda+\mu}c=a\oo{\lambda}(b\oo\mu c), \quad \forall\,
a, b, c\in A.
\end{eqnarray*}
\end{defi}

Let $(A,\; \cdot\ooo{\lambda}{A}\cdot)$, $(B,\; \cdot\ooo{\lambda}{B}\cdot)$ be two
associative conformal algebras. A $\cb[\partial]$-module homomorphism
$f: A\rightarrow B$ is called a {\it homomorphism of associative conformal algebras}
if for any $a, b\in A$, $f(a\ooo{\lambda}{A}b)=f(a)\ooo{\lambda}{B}f(b)$.

\begin{ex}
Let $(A,\cdot)$ be an associative algebra. Then ${\rm Cur}(A)=\cb[\partial]\otimes A$
is an associative conformal algebra with the following $\lambda$-product:
\begin{eqnarray*}
(p(\partial)a)\oo\lambda (q(\partial)b)=p(-\lambda)q(\lambda+\partial)
(a\cdot b), \qquad\forall\, \text{$p(\partial)$, $q(\partial)\in
\cb[\partial]$, $a$, $b\in A$.}
\end{eqnarray*}
\end{ex}

Now we recall that definition of left (or right) modules over an associative
conformal algebra.

\begin{defi}\label{def:mod}
A {\rm (conformal) left module} $M$ over an associative conformal algebra $A$ is a
$\cb[\partial]$-module endowed with a $\cb$-bilinear map $A\times M\rightarrow
M[\lambda]$, $(a, v)\mapsto a\ool{\lambda} v$, satisfying the following axioms:
\begin{align*}
(\partial a)\ool{\lambda} v&=-\lambda a\ool{\lambda} v,\qquad
a\ool{\lambda}(\partial v)=(\partial+\lambda)(a\ool{\lambda} v),\\
&(a\oo{\lambda} b)\ool{\lambda+\mu}v=a\ool{\lambda}(b\ool{\mu} v),
\end{align*}
for any $a, b\in A$ and $v\in M$. We denote it by $(M, \triangleright)$.

A {\rm (conformal) right module} $M$ over an associative conformal algebra $A$ is a
$\cb[\partial]$-module endowed with a $\cb$-bilinear map $M\times A\rightarrow
M[\lambda]$, $(v, a)\mapsto v\oor{\lambda} a$, satisfying
\begin{align*}
(\partial v)\oor{\lambda} a&=-\lambda v\oor{\lambda} a,\qquad
v\oor{\lambda}(\partial a)=(\partial+\lambda)(v\oor{\lambda} a),\\
&(v\oor{\lambda} a)\oor{\lambda+\mu}b=v\oor{\lambda}(a\oo{\mu} b),
\end{align*}
for any $a, b\in A$ and $v\in M$. We denote it by $(M, \triangleleft)$.

An {\rm (conformal) $A$-bimodule} is a triple $(M, \triangleright, \triangleleft)$
such that $(M, \triangleright)$ is a left $A$-module, $(M, \triangleleft)$ is a
right $A$-module, and they satisfy
\begin{eqnarray*}
(a\ool{\lambda} v)\oor{\lambda+\mu}b=a\ool{\lambda}(v\oor{\mu}b),
\end{eqnarray*}
for any $a, b\in A$ and $v\in M$.
\end{defi}

Let $A$ be an associative conformal algebra. Define two linear maps
$\triangleright^{A},\; \triangleleft^{A}:\; A\otimes A\rightarrow A$ by
$a\oll{\lambda}{A} b=a\oo{\lambda} b$ and
$b\orr{\lambda}{A} a=b\oo{\lambda}a$ for all $a$, $b\in A$.
Then $(A, \triangleright^{A}, \triangleleft^{A})$ is a bimodule of $A$,
it is called the regular bimodule of $A$.

In 1999, Bakalov, Kac and Voronov firstly gave the definition of Hochschild
cohomology of associative conformal algebras \cite{BKV}. This definition is a
conformal analogue of Hochschild cohomology of associative algebras.
In \cite{DK}, for the case of Lie conformal algebras, the definition was improved by
taking $n-1$ variables. Following this idea, we define the Hochschild cohomology for
an associative conformal algebra $A$ and a bimodule $M$ over $A$. We denote
$\mathcal{C}^{0}(A, M)=M/\partial M$, $\mathcal{C}^{1}(A, M)=\Hom_{\cb[\partial]}
(A, M)$, the set of $\cb[\partial]$-module homomorphisms from $A$ to $M$, and for
$n\geq2$, the space of $n$-cochains $\mathcal{C}^{n}(A, M)$ consists of all maps
$$
\varphi:\;  A^{\otimes n}\longrightarrow
M[\lambda_{1},\dots, \lambda_{n-1}],
$$
such that
$$
\varphi_{\lambda_{1},\dots, \lambda_{n-1}}(a_{1},\dots, \partial a_{i},\dots, a_{n})
=-\lambda_{i}\varphi_{\lambda_{1},\dots, \lambda_{n-1}}(a_{1},\dots, a_{n}),
$$
for $i=1,2,\cdots,n-1$, and
$$
\varphi_{\lambda_{1},\dots, \lambda_{n-1}}(a_{1},\dots,  a_{n-1},\partial a_{n})
=(\partial+\lambda_{1}+\dots+\lambda_{n-1})\varphi_{\lambda_{1},\dots, \lambda_{n-1}}
(a_{1},\dots, a_{n}).
$$
Such a linear map $\varphi$ is called a {\it conformal sesquilinear map} from
$A^{\otimes n}$ to $M[\lambda_{1},\dots, \lambda_{n-1}]$.
The differentials are defined by $d_{0}: M/\partial M\rightarrow\Hom_{\cb[\partial]}(A, M)$,
\begin{equation*}
d_{0}(v+\partial M)(a)=(a\ool{-\lambda-\partial}v-v\oor{\lambda}a)\mid_{\lambda=0},
\end{equation*}
and for $n\geq1$,
\begin{multline}\nonumber
\qquad d_{n}(\varphi)_{\lambda_{1},\dots, \lambda_{n}}(a_{1},\dots, a_{n+1})
=a_{1}\ool{\lambda_{1}}\varphi_{\lambda_{2},\dots, \lambda_{n}}
(a_{2},\dots, a_{n+1}) \\
+\sum_{i=1}^{n}(-1)^{i}\varphi_{\lambda_{1},\dots, \lambda_{i-1},
\lambda_{i}+\lambda_{i+1}, \lambda_{i+2}\dots,\lambda_{n}}
(a_{1},\dots, a_{i-1}, a_{i}\oo{\lambda_{i}} a_{i+1} , a_{i+2}\dots, a_{n+1}) \\
+(-1)^{n+1}\varphi_{\lambda_{1},\dots, \lambda_{n-1}}(a_{1},\dots,
a_{n})\oor{\lambda_{1}+\dots+\lambda_{n}} a_{n+1}.\qquad
\end{multline}
One can verify that the operator $d_{n}$ preserves the space of cochains and
$d_{n+1}\circ d_{n}=0$. The cochains of an associative conformal algebra $A$ with
coefficients in a bimodule $M$ form a complex $(\mathcal{C}^{\ast}(A, M),
d_{\ast})$, called the {\it Hochschild complex}. We denote the space of
$n$-cocycles by $\mathcal{Z}^{n}(A, M)=\{\varphi\in\mathcal{C}^{n}(A, M)\mid
d_{n}(\varphi)=0\}$, and the space of $n$-coboundaries by $\mathcal{B}^{n}
(A, M)=\{d_{n-1}(\varphi) \mid\varphi\in\mathcal{C}^{n-1}(A, M)\}$. The $n$-th
Hochschild cohomology of $A$ with coefficients in $M$ is define by
\begin{equation*}
\HH^{n}(A, M)=\mathcal{Z}^{n}(A, M)/\mathcal{B}^{n}(A, M).
\end{equation*}
In particular, if $M=A$ as conformal bimodule, we denote $\HH^{n}(A):=\HH^{n}(A, M)$,
and call the $n$-th Hochschild cohomology of $A$. In this paper, we consider
the Gerstenhaber algebra structure on $\HH^{n}(A)$.
For the description of the lowest degree of the Hochschild cohomology, we have
\begin{itemize}
\item[(1)] $\Img d_{0}=\Inn(A, M)$, where $\Inn(A, M)=\{f_{v}\in\Hom_{\cb[\partial]}
           (A, M)\mid v\in M, f_{v}(a)=a\oo{-\partial}v-v\oo{0}a\}$;
\item[(2)] $\Ker d_{1}=\Der(A, M)$, where $\Der(A, M)=\{f\in\Hom_{\cb[\partial]}
           (A, M)\mid f(a\oo{\lambda}b)=a\ool{\lambda}f(b)+f(a)\oor{\lambda}b\}$;
\item[(3)] if $A$ as $\cb[\partial]$-module is projective, the equivalence classes
           of $\cb[\partial]$-split abelian extensions of $A$ by conformal $A$-bimodule $M$
           are correspond bijectively to $\HH^{n}(A, M)$. For the details see \cite{Do,Lib}.
\end{itemize}


\section{Graded Lie algebraic structure}\label{sec:lie}

In this section, we give a graded Lie algebraic structure on the Hochschild cohomology
of an associative conformal algebra. First, let us recall the Gerstenhaber's
classic conclusion for the pre-Lie system.

\begin{defi}\label{def:pre-Lie sys}
A pre-Lie system $\{V_{m}, \circ_{i}\}$ consists of a sequence of vector spaces
$V_{m}$, $m \geq 0$, and for each $m , n \geq 0$ there exist linear maps
\begin{equation*}
\circ_{i} : V_{m} \otimes V_{n} \rightarrow V_{m+n}, \qquad
(f, g)\mapsto f \circ_{i} g,
\end{equation*}
for any $0\leq i \leq m$, satisfying
\begin{equation*}
(f\circ_{i} g)\circ_{j} h=
\begin{cases}(f \circ_{j} h) \circ_{i+p} g, &\mbox{ if} \quad 0\leq j\leq i-1;\\
f\circ_{i} (g \circ_{j-i} h), &\mbox{ if} \quad i\leq j\leq n+i,
\end{cases}
\end{equation*}
for $f\in V_{m}$, $g\in V_{n}$ and $h\in V_{p}$.
\end{defi}

Let $\{V_{m}, \circ_{i}\}$ be a pre-Lie system. Then for any $m, n \geq 0$,
there is a new linear map $\circ : V_m \otimes V_n \rightarrow V_{m+n}$,
$(f, g)\mapsto f\bar{\circ} g$ defined by
\begin{equation*}
f\bar{\circ} g = \sum_{i=0}^{m}(-1)^{ni}f \circ_{i} g,
\end{equation*}
for any $f\in V_{m}$, $g \in V_{n}$. Then we have the following theorem.

\begin{thm}[\cite{Ger}]\label{thm: Pre-Lie alg}
Let $\{V_{m}, \circ_{i}\}$ be a pre-Lie system and $f\in V_{m}$, $g \in V_{n}$
and $h \in V_{p}$, respectively. Then
\begin{itemize}
\item[(i)] $(f\bar{\circ} g)\bar{\circ} h-f\bar{\circ}(g\bar{\circ} h)
           =\sum_{i, j}(-1)^{ni+pj}(f\circ_{i}g)\circ_{j}h$, where the summation
           is indexed over those $i$ and $j$ with either $0\leq j\leq i-1$ or
           $n+i+1\leq j\leq m+n$,
\item[(ii)] $(f\bar{\circ} g)\bar{\circ} h-f\bar{\circ}(g\bar{\circ} h)
           =(-1)^{np}[(f\bar{\circ} h)\bar{\circ} g-f\bar{\circ}(h\bar{\circ} g)]$.
\end{itemize}
\end{thm}

\begin{defi}\label{def:pre-Lie alg}
A {\rm graded pre-Lie algebra} $\{V_{m}, \bar{\circ}\}$ consists of vector spaces
$V_{m}$ together with linear maps $\bar{\circ}: V_{m}\otimes V_{n}\rightarrow V_{m+n}$
satisfying
\begin{equation*}
(f\bar{\circ} g)\bar{\circ} h-f\bar{\circ}(g\bar{\circ} h)
=(-1)^{np}[(f\bar{\circ} h)\bar{\circ} g-f\bar{\circ}(h\bar{\circ} g)],
\end{equation*}
for any $f\in V_{m}$, $g\in V_{n}$ and $h\in V_{p}$.
\end{defi}

It is proved in \cite{Ger} that the graded commutator of a graded pre-Lie algebra
defines a graded Lie algebra structure. More precisely, if $\{V_{m}, \bar{\circ}\}$
is a graded pre-Lie algebra, then the bracket $[-,-]$ defined by
$[f, g]=f\bar{\circ} g-(-1)^{mn}g\bar{\circ} f$,  for any $f\in V_{m}$, $g \in V_{n}$,
defines a graded Lie algebra structure on $\bigoplus_{m\geq0}V_{m}$.

Let now $A$ be an associative conformal algebra.
We denote by $U_{0}=\Hom_{\cb[\partial]}(A, A)$, and denote by $U_{n}$ the set of
all conformal sesquilinear maps from $A^{\otimes n+1}$ to $A[\lambda_{0},\cdots,
\lambda_{n-1}]$ for $n\geq1$. Define $\bullet_{i}: U_{m}\times U_{n}\rightarrow U_{m+n}$
by
\begin{align*}
&\;(f\bullet_{i} g)_{\lambda_{0},\dots,\lambda_{m+n-1}}(a_{0}, a_{1},\dots, a_{m+n})\\
=&\;f_{\lambda_{0},\dots,\lambda_{i-1}, \lambda_{i}+\cdots+\lambda_{i+n},
\lambda_{i+n+1},\dots,\lambda_{m+n-1}}\big(a_{0},\dots, a_{i-1}, g_{\lambda_{i},
\dots,\lambda_{i+n-1}}(a_{i},\dots, a_{i+n}), a_{i+n+1}, \dots, a_{m+n}\big),
\end{align*}
for any $f\in U_{m}$, $g\in U_{n}$ and $0\leq i\leq m$.
In particular, if $m=0$,
\begin{equation*}
(f\bullet_{0} g)_{\lambda_{0},\dots,\lambda_{n-1}}(a_{0}, a_{1},
\dots, a_{n})=f\Big(g_{\lambda_{0},\dots,\lambda_{n-1}}(a_{0}, a_{1},
\dots, a_{n})\Big),
\end{equation*}
where $f$  is extended canonically to a $\cb[\partial]$-module homomorphism from
$A[\lambda_{0},\dots,\lambda_{n-1}]$ to $A[\lambda_{0},\dots,\lambda_{n-1}]$,
and if $n=0$,
\begin{equation*}
(f\bullet_{i} g)_{\lambda_{0},\dots,\lambda_{m-1}}
(a_{0}, a_{1},\dots, a_{m})=f_{\lambda_{0},\dots,\lambda_{m-1}}
(a_{0},\dots, a_{i-1}, g(a_{i}), a_{i-1},\dots, a_{m}).
\end{equation*}
Then we get a pre-Lie system.

\begin{lem}\label{lem:pre-Lie sys}
With the above notations, $\{U_{m}, \bullet_{i}\}$ forms a pre-Lie system.
\end{lem}

\begin{proof}
Direct calculation shows that for any $f\in U_{m}$, $g\in U_{n}$ and $h\in U_{p}$,
(1) if $0\leq j\leq i-1$, we have
\begin{align*}
&\;((f\bullet_{j}h)\bullet_{i+p}g)_{\lambda_{0},\dots,\lambda_{m+n+p-1}}
(a_{0}, a_{1},\dots, a_{m+n+p})\\
=&\;(f\bullet_{j}h)_{\lambda_{0},\dots,\lambda_{i+p-1}, \lambda_{i+p}+\cdots
+\lambda_{i+p+n},\lambda_{i+p+n+1},\dots,\lambda_{m+n+p-1}}
\big(a_{0},\dots, a_{i+p-1}, \\
&\qquad g_{\lambda_{i+p},\dots,\lambda_{i+p+n-1}}(a_{i+p},\dots, a_{i+n+p}),
a_{i+n+p+1}, \dots, a_{m+n+p}\big)\\
=&\;f_{\lambda_{0},\dots,\lambda_{j-1}, \lambda_{j}+\cdots+\lambda_{j+p},
\lambda_{j+p+1},\dots,\lambda_{i+p-1},\lambda_{i+p}+\cdots
+\lambda_{i+p+n},\lambda_{i+p+n+1},\dots,\lambda_{m+n+p-1}}
\big(a_{0},\dots, a_{j-1}, h_{\lambda_{j},\dots,\lambda_{j+p-1}}
(a_{j},\dots,a_{j+p})\\
&\qquad a_{j+p+1},\dots, a_{i+p-1}, g_{\lambda_{i+p},\dots,\lambda_{i+p+n-1}}
(a_{i+p},\dots, a_{i+n+p}),a_{i+n+p+1}, \dots, a_{m+n+p}\big),\\
=&\;((f\bullet_{i} g)\bullet_{j} h)_{\lambda_{0},\dots,\lambda_{m+n+p-1}}
(a_{0}, a_{1},\dots, a_{m+n+p});
\end{align*}
(2) if $i\leq j\leq n+i$, we have
\begin{align*}
&\;(f\bullet_{i}(g\bullet_{j-i}h))_{\lambda_{0},\dots,\lambda_{m+n+p-1}}
(a_{0}, a_{1},\dots, a_{m+n+p})\\
=&\;f_{\lambda_{0},\dots,\lambda_{i-1},\lambda_{i}+\cdots
+\lambda_{i+p+n},\lambda_{i+p+n+1},\dots,\lambda_{m+n+p-1}}
\big(a_{0},\dots, a_{i-1}, \\
&\qquad (g\bullet_{j-i}h)_{\lambda_{i},\dots,\lambda_{i+p+n-1}}(a_{i},
\dots, a_{i+n+p}),a_{i+n+p+1}, \dots, a_{m+n+p}\big)\\
=&\;f_{\lambda_{0},\dots,\lambda_{i-1},\lambda_{i}+\cdots
+\lambda_{i+p+n},\lambda_{i+p+n+1},\dots,\lambda_{m+n+p-1}}
\big(a_{0},\dots, a_{i-1}, g_{\lambda_{i},\dots,\lambda_{j-1}, \lambda_{j}+\cdots+
\lambda_{j+p}, \lambda_{j+p+1},\dots, \lambda_{i+n+p-1}}(a_{i},\dots,\\
&\qquad a_{j-1}, h_{\lambda_{j},\dots,\lambda_{j+p-1}}(a_{j},\dots,a_{j+p}),
a_{j+p+1},\dots, a_{i+n+p}), a_{i+n+p+1}, \dots, a_{m+n+p}\big)\\
=&\;((f\bullet_{i} g)\bullet_{j} h)_{\lambda_{0},\dots,\lambda_{m+n+p-1}}
(a_{0}, a_{1},\dots, a_{m+n+p}).
\end{align*}
Thus $\{U_{m}, \bullet_{i}\}$ is a pre-Lie system.
\end{proof}

Then by Theorem \ref{thm: Pre-Lie alg}, we get a pre-Lie algebra
$\Big(\bigoplus_{i\geq0}U_{i},\; \bullet\Big)$, where
\begin{equation*}
f\bullet g = \sum_{i=0}^{m}(-1)^{ni}f\bullet_{i} g,
\end{equation*}
for any $f\in U_{m}$, $g\in U_{n}$. Therefore, we get a graded Lie algebra
$\Big(\bigoplus_{i\geq0}U_{i},\; [-, -]\Big)$ with
\begin{equation*}
[f, g]=f\bullet g-(-1)^{mn}g\bullet f,\qquad \forall\; f\in U_{m},\; g\in U_{n}.
\end{equation*}
One can check that $\rho\in U_{2}$ such that $(A, \rho)$ is an associative
conformal algebra if and only if $[\rho, \rho]=0$, that is, $\rho$ satisfying
the Maurer-Cartan equation of $\Big(\bigoplus_{i\geq0}U_{i},\; [-, -]\Big)$.

Let $A$ be an associative conformal algebra, where we denote the conformal
multiplication by $\rho(a, b)=a\oo{\lambda}b$ for any $a, b\in A$.
Now go back to the Hochschild cohomology of associative conformal algebra $A$.
For any $b\in\mathcal{C}^{0}(A, A)$ and $f\in\mathcal{C}^{m}(A, A)$, we define
$b\bullet f=0$ and $f\bullet b$ by
\begin{align*}
(f\bullet b)_{\lambda_{1},\dots,\lambda_{m-2}}(a_{1},\dots,a_{m-1})
&=\sum_{i=1}^{m-1}(-1)^{i-1}f_{\lambda_{1},\dots,\lambda_{i-1}, 0, \lambda_{i},\dots,
\lambda_{m-2}}(a_{1},\dots,a_{i-1}, b, a_{i},\dots,a_{m-1})\\
&\qquad\qquad +(-1)^{m-1}f_{\lambda_{1},\dots,\lambda_{m-2}, -\partial}(a_{1},\dots,a_{m-1}, b).
\end{align*}
We define degree $-1$ bracket $[-, -]: \mathcal{C}^{m}(A, A)\times
\mathcal{C}^{n}(A, A)\rightarrow\mathcal{C}^{m+n-1}(A, A)$ by
\begin{equation*}
[f, g]=f\bullet g-(-1)^{(m-1)(n-1)}g\bullet f,
\end{equation*}
for any $f\in\mathcal{C}^{m}(A, A)$ and $g\in\mathcal{C}^{n}(A, A)$, $m, n\geq0$.
First, note the space of $n$-cochains $\mathcal{C}^{n}(A, A)=U_{n-1}$ for all
$n\geq1$, we get a graded Lie algebra $\Big(\bigoplus_{i\geq1}\mathcal{C}^{i}(A, A),\;
[-, -]\Big)$. Second, for any $b, b'\in\mathcal{C}^{0}(A, A)$, $f\in\mathcal{C}^{m}(A, A)$,
$g\in\mathcal{C}^{n}(A, A)$ and $h\in\mathcal{C}^{p}(A, A)$,
one can check that $[b, b']=0$, $[f, b]=(-1)^{m-1}[b, f]$, and
$[f, [g, h]]=[[f, g], h]+(-1)^{(m-1)(n-1)}[g, [f, h]]$ if at least one
of $m, n$ and $p$ is zero, by direct calculation. Thus,
$\bigoplus_{i\geq0}\mathcal{C}^{i}(A, A)$ is a graded Lie algebra under the
degree $-1$ bracket $[-, -]$. For any cocycle $f\in\mathcal{C}^{m}(A, A)$,
we equate $f$ with its image in $\HH^{n}(A)$, then the Lie bracket induces a
multiplication $[-, -]:\HH^{m}(A)\times \HH^{m}(A)\rightarrow \HH^{m+n-1}(A)$.

\begin{thm}\label{thm:graded-Lie}
The Lie bracket $[-, -]$ on $\bigoplus_{i\geq0}\mathcal{C}^{i}(A, A)$
induces a bracket $[-, -]$ of degree $-1$ on $\bigoplus_{i\geq0}\HH^{i}(A)$ such that
$\Big(\bigoplus_{i\geq0}\HH^{i}(A),\; [-, -]\Big)$ is a graded Lie algebra.
\end{thm}

\begin{proof}
We denote the conformal multiplication on $A$ by $\rho$, i.e.,
$\rho_{\lambda}(a, b)=a\oo{\lambda}b$. Then for any $f\in\mathcal{C}^{m}(A, A)$,
we have
\begin{align*}
&\;(f\bullet\rho-(-1)^{m-1}\rho\bullet f)_{\lambda_{1},\dots,\lambda_{m}}
(a_{1},\dots, a_{m+1})\\
=&\; \sum_{i=1}^{m}(-1)^{i-1}f_{\lambda_{1},\dots,\lambda_{i-1}, \lambda_{i}+\lambda_{i+1},
\lambda_{i+2}\dots,\lambda_{m}}\Big(a_{1},\dots,a_{i-1},\rho_{\lambda_{i}}(a_{i},
a_{i+1}),a_{i+2},\dots,a_{m+1}\Big)\\
&\;-(-1)^{m-1}\rho_{\lambda_{1}+\cdots+\lambda_{m}}\Big(f_{\lambda_{1},\dots,
\lambda_{m-1}}(a_{1},\dots,a_{m}), a_{m+1}\Big)
-\rho_{\lambda_{1}}\Big(a_{1}, f_{\lambda_{2},\dots,\lambda_{m}}(a_{2},
\dots,a_{m+1})\Big)\\
=&\;-d_{m}(f)_{\lambda_{1},\dots,\lambda_{m}}(a_{1},\dots, a_{m+1}).
\end{align*}
That is, $d_{m}(f)=-[f, \rho]=(-1)^{m+1}[\rho, f]$. Thus, for any
$f\in\mathcal{C}^{m}(A, A)$ and $g\in\mathcal{C}^{n}(A, A)$, $m, n\geq1$,
\begin{align*}
d_{m+n-1}([f, g])=&(-1)^{m+n}[\rho, [f, g]]\\
=&(-1)^{m+n}\Big([[\rho, f], g]]+(-1)^{n+1}[f, [\rho, g]]\Big)\\
=&(-1)^{n+1}[d_{m}(f), g]+[f, d_{n}(g)].
\end{align*}
Moreover, by direct calculation, we also have $d_{m-1}([f, b])=[f, d_{0}(b)]-[d_{m}(f), b]$
for $b\in\mathcal{C}^{0}(A, A)$ and $f\in\mathcal{C}^{m}(A, A)$.
This means that if $f\in\mathcal{Z}^{m}(A, A)$ and $g\in\mathcal{Z}^{n}(A, A)$,
then $[f, g]\in\mathcal{Z}^{m+n-1}(A, A)$; if $f\in\mathcal{Z}^{m}(A, A)$,
$g\in\mathcal{Z}^{n}(A, A)$, and $f\in\mathcal{B}^{m}(A, A)$ or
$g\in\mathcal{B}^{n}(A, A)$, then $[f, g]\in\mathcal{B}^{m+n-1}(A, A)$, for any
$m, n\geq0$. Hence the bracket $[-, -]$ is well-defined,
and so that $\Big(\bigoplus_{i\geq0}\HH^{i}(A),\; [-, -]\Big)$ is a graded Lie algebra.
\end{proof}

Under the bracket $[-, -]$, $\HH^{1}(A)$ is a Lie algebra, and $\HH^{n}(A)$
are $\HH^{1}(A)$-modules for all $n\geq2$.


\section{The cup product}\label{sec:cup}

In this section, we define a product on the Hochschild cohomology complex
of an associative conformal algebra. This product is conformal analogue of
the cup product on Hochschild cohomology complex of associative algebras.
Thus we also call it cup product. For any $f\in\mathcal{C}^{m}(A, A)$ and
$g\in\mathcal{C}^{n}(A, A)$, $m, n\geq1$, we define $f\sqcup g\in\mathcal{C}^{m+n}(A, A)$ by
\begin{equation*}
(f\sqcup g)_{\lambda_{1},\dots,\lambda_{m+n-1}}(a_{1},\dots, a_{m+n})
=f_{\lambda_{1},\dots,\lambda_{m-1}}(a_{1},\dots, a_{m})\oo{\lambda_{1}+\cdots+\lambda_{m}}
g_{\lambda_{m+1},\dots,\lambda_{m+n-1}}(a_{m+1},\dots, a_{m+n}).
\end{equation*}
If we view the conformal multiplication $\rho$ on $A$ as an element
in $\mathcal{C}^{2}(A, A)$, we have
\begin{equation*}
f\sqcup g=((\rho\bullet_{0}f)\bullet_{m} g),\qquad
g\sqcup f=((\rho\bullet_{1}f)\bullet_{0} g).
\end{equation*}
If $b\in\mathcal{C}^{0}(A, A)$, we define
\begin{align*}
(b\sqcup g)_{\lambda_{1},\dots,\lambda_{n-1}}(a_{1},\dots, a_{n})
&=b\oo{0}g_{\lambda_{1},\dots,\lambda_{n-1}}(a_{1},\dots, a_{n}),\\
(f\sqcup b)_{\lambda_{1},\dots,\lambda_{m-1}}(a_{1},\dots, a_{m})
&=f_{\lambda_{1},\dots,\lambda_{m-1}}(a_{1},\dots, a_{m})\oo{-\partial}b.
\end{align*}
Then for any $f\in\mathcal{C}^{m}(A, A)$, $g\in\mathcal{C}^{n}(A, A)$ and
$h\in\mathcal{C}^{p}(A, A)$, $m, n, p\geq1$, we have
\begin{align*}
&\;(f\sqcup(g\sqcup h))_{\lambda_{1},\dots,\lambda_{m+n+p-1}}(a_{1},\dots, a_{m+n+p})\\
=&\;f_{\lambda_{1},\dots,\lambda_{m-1}}(a_{1},\dots, a_{m})\oo{\lambda_{1}
+\cdots+\lambda_{m}}(g\sqcup h)_{\lambda_{m+1},\dots,\lambda_{m+n+p-1}}
(a_{m+1},\dots, a_{m+n+p})\\
=&\;f_{\lambda_{1},\dots,\lambda_{m-1}}(a_{1},\dots, a_{m})\oo{\lambda_{1}
+\cdots+\lambda_{m}}\Big(g_{\lambda_{m+1},\dots,\lambda_{m+n-1}}
(a_{m+1},\dots, a_{m+n})\oo{\lambda_{m+1}+\cdots+\lambda_{m+n}}\\
&\qquad h_{\lambda_{m+n+1},\dots,\lambda_{m+n+p-1}}(a_{m+n+1},\dots, a_{m+n+p})\Big)\\
=&\;\Big(f_{\lambda_{1},\dots,\lambda_{m-1}}(a_{1},\dots, a_{m})\oo{\lambda_{1}
+\cdots+\lambda_{m}}g_{\lambda_{m+1},\dots,\lambda_{m+n-1}}(a_{m+1},\dots,
a_{m+n})\Big)\oo{\lambda_{1}+\cdots+\lambda_{m+n}}\\
&\qquad h_{\lambda_{m+n+1},\dots,\lambda_{m+n+p-1}}(a_{m+n+1},\dots, a_{m+n+p})\\
=&\;((f\sqcup g)\sqcup h)_{\lambda_{1},\dots,\lambda_{m+n+p-1}}(a_{1},\dots, a_{m+n+p}).
\end{align*}
Similarly, one can check that $f\sqcup(g\sqcup h)=(f\sqcup g)\sqcup h$ for
$f\in\mathcal{C}^{m}(A, A)$, $g\in\mathcal{C}^{n}(A, A)$ and
$h\in\mathcal{C}^{p}(A, A)$, $m, n, p\geq0$.
That is, $\Big(\bigoplus_{i\geq0}\mathcal{C}^{i}(A, A),\; \sqcup\Big)$ is a
graded associative algebra. We call the product $\sqcup$ the cup product.
The following proposition show that the cup product is also inducing a
graded associative algebra structure on $\bigoplus_{i\geq1}\HH^{i}(A)$.

\begin{pro}\label{pro: cup-product}
The cup product $\sqcup$ on $\bigoplus_{i\geq0}\mathcal{C}^{i}(A, A)$ induces a
product on $\bigoplus_{i\geq0}\HH^{i}(A)$, denoted by $\sqcup$, such that
$\Big(\bigoplus_{i\geq0}\HH^{i}(A),\; \sqcup\Big)$ is a graded associative algebra.
\end{pro}

\begin{proof}
We firstly show that the differential $d$ on complex $\mathcal{C}^{\bullet}(A, A)$
satisfies the graded Leibniz rule with respect to the cup product.
For any $f\in\mathcal{C}^{m}(A, A)$, $g\in\mathcal{C}^{n}(A, A)$, $m, n\geq1$,
we have
\begin{align*}
&\;d_{m+n}(f\sqcup g)_{\lambda_{1},\dots,\lambda_{m+n}}(a_{1},\dots, a_{m+n+1})\\
=&\;a_{1}\oo{\lambda_{1}}(f\sqcup g)_{\lambda_{2},\dots,\lambda_{m+n}}
(a_{2},\dots, a_{m+n+1})\\
&\;+\sum_{i=1}^{m+n}(-1)^{i}(f\sqcup g)_{\lambda_{1},\dots,\lambda_{i-1},
\lambda_{i}+\lambda_{i+1}, \lambda_{i+2},\dots,\lambda_{m+n}}(a_{1},
\dots, a_{i-1}, a_{i}\oo{\lambda_{i}}a_{i+1}, a_{i+2},\dots, a_{m+n+1})\\
&\;+(-1)^{m+n+1}(f\sqcup g)_{\lambda_{1},\dots,\lambda_{m+n-1}}
(a_{1},\dots, a_{m+n})\oo{\lambda_{1}+\cdots+\lambda_{m+n}}a_{m+n+1}\\
=&\;a_{1}\oo{\lambda_{1}}\Big(f_{\lambda_{2},\dots,\lambda_{m}}(a_{2},\dots, a_{m+1})
\oo{\lambda_{2}+\cdots+\lambda_{m+1}}g_{\lambda_{m+2},\dots,\lambda_{m+n}}
(a_{m+2},\dots, a_{m+n+1})\Big)\\
&\;+\sum_{i=1}^{m}(-1)^{i}f_{\lambda_{1},\dots,\lambda_{i-1},\lambda_{i}
+\lambda_{i+1}, \lambda_{i+2},\cdots,\lambda_{m}}(a_{1}, \dots, a_{i-1}, a_{i}
\oo{\lambda_{i}}a_{i+1}, a_{i+2},\dots, a_{m+1})\oo{\lambda_{1}
+\cdots+\lambda_{m+1}}\qquad\qquad\\
&\qquad g_{\lambda_{m+2},\cdots,\lambda_{m+n}}(a_{m+2},\dots, a_{m+n+1})
\end{align*}
\begin{align*}
&\;+\sum_{i=m+1}^{m+n}(-1)^{i}f_{\lambda_{1},\dots,\lambda_{m-1}}(a_{1},\dots,
a_{m})\oo{\lambda_{1}+\cdots+\lambda_{m}}g_{\lambda_{m+1},\dots,\lambda_{i-1},
\lambda_{i}+\lambda_{i+1}, \lambda_{i+2},\dots,\lambda_{m+n}}(a_{m+1}, \dots,\\
&\qquad a_{i-1}, a_{i}\oo{\lambda_{i}}a_{i+1}, a_{i+2},\dots, a_{m+n+1})\\
&\;+(-1)^{m+n+1}\Big(f_{\lambda_{1},\dots,\lambda_{m-1}}(a_{1},\dots, a_{m})
\oo{\lambda_{1}+\cdots+\lambda_{m}}g_{\lambda_{m+1},\dots,\lambda_{m+n-1}}
(a_{m+1},\dots, a_{m+n})\Big)\oo{\lambda_{1}+\cdots+\lambda_{m+n}}a_{m+n+1}\\
=&\;d_{m}(f)_{\lambda_{1},\dots,\lambda_{m}}(a_{1},\dots, a_{m+1})\oo{
\lambda_{1}+\cdots+\lambda_{m+1}}g_{\lambda_{m+2},\dots,\lambda_{m+n}}
(a_{m+2},\dots, a_{m+n+1})\\
&\qquad +(-1)^{m}f_{\lambda_{1},\dots,\lambda_{m-1}}(a_{1},\dots, a_{m})\oo{
\lambda_{1}+\cdots+\lambda_{m}}d_{n}(g)_{\lambda_{m+1},\dots,\lambda_{m+n}}
(a_{m+1},\dots, a_{m+n+1})\\
=&\;\Big(d_{m}(f)\sqcup g+(-1)^{m}f\sqcup d_{n}(g)\Big)_{\lambda_{1},\dots,\lambda_{m+n}}
(a_{1},\dots, a_{m+n+1}).
\end{align*}
Thus $d_{m+n}(f\sqcup g)=d_{m}(f)\sqcup g+(-1)^{m}f\sqcup d_{n}(g)$.
Similarly, one can get the equation $d_{m+n}(f\sqcup g)=d_{m}(f)\sqcup g
+(-1)^{m}f\sqcup d_{n}(g)$ is also true whenever $m=0$ or $n=0$ by direct calculation.
This means that if $f\in\mathcal{Z}^{m}(A, A)$ and $g\in\mathcal{Z}^{n}(A, A)$,
then $f\sqcup g\in\mathcal{Z}^{m+n}(A, A)$; if $f\in\mathcal{Z}^{m}(A, A)$,
$g\in\mathcal{Z}^{n}(A, A)$, and $f\in\mathcal{B}^{m}(A, A)$ or
$g\in\mathcal{B}^{n}(A, A)$, then $f\sqcup g\in\mathcal{B}^{m+n}(A, A)$.
Hence the cup product $\sqcup$ is well-defined,
and $\Big(\bigoplus_{i\geq1}\HH^{i}(A),\; \sqcup\Big)$ is a graded
associative algebra.
\end{proof}

In general, the cup product $\sqcup$ on the Hochschild cohomology complex
$\mathcal{C}^{\ast}(A, A)$ is not graded commutative. But we have the following
theorem.

\begin{thm}\label{thm: commutative}
The graded associative algebra $\Big(\bigoplus_{i\geq0}\HH^{i}(A),\; \sqcup\Big)$
is graded commutative.
\end{thm}

\begin{proof}
For any $f\in\mathcal{C}^{m}(A, A)$, $g\in\mathcal{C}^{n}(A, A)$, $m, n\geq1$,
we have a equation
\begin{equation*}
f\bullet d_{n}(g)+(-1)^{n-1}d_{m}(f)\bullet g-d_{m+n-1}(f\bullet g)
=(-1)^{n-1}\Big(g\sqcup f-(-1)^{mn}f\sqcup g\Big).
\end{equation*}
Indeed, note that
\begin{align*}
&\;f\bullet d_{n}(g)+(-1)^{n-1}d_{m}(f)\bullet g-d_{m+n-1}(f\bullet g)\\
=&\;(-1)^{n-1}f\bullet(\rho\bullet g)-f\bullet(g\bullet\rho)
+(-1)^{m+n}(\rho\bullet f)\bullet g\\
&\qquad -(-1)^{n-1}(f\bullet\rho)\bullet g-(-1)^{m+n}\rho\bullet(f\bullet g)
+(f\bullet g)\bullet\rho,
\end{align*}
and $\Big(\bigoplus_{i\geq0}U_{i},\; \bullet\Big)$ is a pre-Lie algebra, i.e.,
\begin{equation*}
(f\bullet g)\bullet\rho-f\bullet(g\bullet\rho)
=(-1)^{n-1}\Big((f\bullet\rho)\bullet g-f\bullet(\rho\bullet g)\Big),
\end{equation*}
we get
\begin{align*}
&\;f\bullet d_{n}(g)+(-1)^{n-1}d_{m}(f)\bullet g-d_{m+n-1}(f\bullet g) \\
=&\;(-1)^{m+n}\Big((\rho\bullet f)\bullet g-\rho\bullet(f\bullet g)\Big) \\
=&\;(-1)^{m+n}\Big((-1)^{m(n-1)}(\rho\bullet_{0}f)\bullet_{m}g
+(-1)^{m-1}(\rho\bullet_{1}f)\bullet_{0} g\Big) \\
=&\;(-1)^{m+n}\Big((-1)^{m(n-1)}f\sqcup g+(-1)^{m-1}g\sqcup f\Big)\\
=&\;(-1)^{n-1}\Big(g\sqcup f-(-1)^{mn}f\sqcup g\Big).
\end{align*}
Thus, if $f\in\mathcal{Z}^{m}(A, A)$ and $g\in\mathcal{Z}^{n}(A, A)$, then
\begin{equation*}
(-1)^{n-1}\Big(g\sqcup f-(-1)^{mn}f\sqcup g\Big)=d_{m+n-1}(f\bullet g)
\in\mathcal{B}^{m+n-1}(A, A).
\end{equation*}
Hence, $g\sqcup f=(-1)^{mn}f\sqcup g$ if $f\in\HH^{m}(A)$ and $g\in\HH^{n}(A)$.
Note that $b\oo{0}a=a\oo{-\partial}b$ for any $a\in A$ if $b\in\HH^{0}(A)$,
we obtain $g\sqcup f=(-1)^{mn}f\sqcup g$ if $m=0$ or $n=0$.
Thus, $\Big(\bigoplus_{i\geq0}\HH^{i}(A),\; \sqcup\Big)$ is graded commutative.
\end{proof}


\section{The Gerstenhaber algebraic structure}\label{sec:Gerst}

In this section, we prove that the degree $-1$ graded Lie bracket $[-,-]$ and
the cup product $\sqcup$ on the Hochschild cohomology $\bigoplus_{i\geq0}\HH^{i}(A)$
of an associative conformal algebra $A$ satisfy the Leibniz rule of a
Gerstenhaber algebra. This gives the Gerstenhaber algebraic structure on
$\bigoplus_{i\geq0}\HH^{i}(A)$. For any $f\in\mathcal{C}^{m}(A, A)$,
$g\in\mathcal{C}^{n}(A, A)$, $h\in\mathcal{C}^{p}(A, A)$, $m, n, p\geq1$,
$1\leq i\leq p-1$ and $m+i\leq j\leq m+p-1$, we set
\begin{align*}
H_{i,j}=& a_{1}\oo{\lambda_{1}}h_{\lambda_{2},\dots,\lambda_{i},
\lambda_{i+1}+\cdots+\lambda_{i+m},\lambda_{i+m+1},\dots,\lambda_{j},
\lambda_{j+1}+\cdots+\lambda_{j+n},\lambda_{j+n+1},\dots,\lambda_{m+n+p-2}}
\Big(a_{2},\dots, a_{i}, \\
&\quad f_{\lambda_{i+1},\dots,\lambda_{i+m-1}}(a_{i+1},\dots, a_{i+m}),
a_{i+m+1},\dots, a_{j},\\
&\quad g_{\lambda_{j+1},\dots,\lambda_{j+n-1}}(a_{j+1},
\dots, a_{j+n}), a_{j+n+1},\dots, a_{m+n+p-1}\Big)\\
&+\sum_{q=1}^{i-1}(-1)^{q}h_{\lambda_{1},\dots,\lambda_{q-1}, \lambda_{q}
+\lambda_{q+1},\lambda_{q+2},\dots,\lambda_{i}, \lambda_{i+1}+\cdots+\lambda_{i+m},
\lambda_{i+m+1},\dots,\lambda_{j},\lambda_{j+1}+\cdots+\lambda_{j+n},
\lambda_{j+n+1},\dots, \lambda_{m+n+p-2}}\\
&\quad \Big(a_{1},\dots,a_{q-1}, a_{q}\oo{\lambda_{q}}a_{q+1}, a_{q+2},\dots, a_{i},
f_{\lambda_{i+1},\dots,\lambda_{i+m-1}}(a_{i+1},\dots, a_{i+m}),\\
&\quad a_{i+m+1},\dots,a_{j}, g_{\lambda_{j+1},\dots,\lambda_{j+n-1}}
(a_{j+1},\dots, a_{j+n}),a_{j+n+1},\dots, a_{m+n+p-1}\Big)\\
&+(-1)^{i}h_{\lambda_{1},\dots,\lambda_{i-1},
\lambda_{i}+\cdots+\lambda_{i+m},\lambda_{i+m+1},\dots,\lambda_{j},
\lambda_{j+1}+\cdots+\lambda_{j+n},\lambda_{j+n+1},\dots,\lambda_{m+n+p-2}}
\Big(a_{1},\dots, a_{i-1}, \\
&\quad a_{i}\oo{\lambda_{i}}f_{\lambda_{i+1},\dots,\lambda_{i+m-1}}(a_{i+1},\dots,
a_{i+m}), a_{i+m+1},\dots, a_{j},\\
&\quad g_{\lambda_{j+1},\dots,\lambda_{j+n-1}}(a_{j+1},
\dots, a_{j+n}), a_{j+n+1},\dots, a_{m+n+p-1}\Big),
\end{align*}
\begin{align*}
H_{i,j}'=&(-1)^{m+i-1}h_{\lambda_{1},\dots,\lambda_{i-1},
\lambda_{i}+\cdots+\lambda_{i+m},\lambda_{i+m+1},\dots,\lambda_{j},
\lambda_{j+1}+\cdots+\lambda_{j+n},\lambda_{j+n+1},\dots,\lambda_{m+n+p-2}}
\Big(a_{1},\dots, a_{i-1},\\
&\quad f_{\lambda_{i},\dots,\lambda_{i+m-2}}
(a_{i},\dots, a_{i+m-1})\oo{\lambda_{i}+\cdots+\lambda_{i+m-1}}
a_{i+m}, a_{i+m+1},\dots, a_{j}, \\
&\quad g_{\lambda_{j+1},\dots,\lambda_{j+n-1}}(a_{j+1},
\dots, a_{j+n}), a_{j+n+1},\dots, a_{m+n+p-1}\Big) \\
&+\sum_{q=m+i}^{j-1}(-1)^{q}h_{\lambda_{1},\dots,\lambda_{i-1},
\lambda_{i}+\cdots+\lambda_{i+m-1}, \lambda_{i+m},\dots,\lambda_{q-1},
\lambda_{q}+\lambda_{q+1},\lambda_{q+2},\dots,\lambda_{j},\lambda_{j+1}+\cdots
+\lambda_{j+n},\lambda_{j+n+1},\dots,\lambda_{m+n+p-2}}\\
&\quad \Big(a_{1},\dots,a_{i-1}, f_{\lambda_{i},\dots,\lambda_{i+m-2}}
(a_{i},\dots, a_{i+m-1}), a_{i+m},\dots, a_{q-1}, a_{q}\oo{\lambda_{q}}a_{q+1},\\
&\quad a_{q+2},\dots, a_{j}, g_{\lambda_{j+1},\dots,\lambda_{j+n-1}}(a_{j+1},
\dots, a_{j+n}), a_{j+n+1},\dots, a_{m+n+p-1}\Big)\\
&+(-1)^{j}h_{\lambda_{1},\dots,\lambda_{i-1},
\lambda_{i}+\cdots+\lambda_{i+m-1}, \lambda_{i+m},\dots,\lambda_{j-1},
\lambda_{j}+\cdots+\lambda_{j+n},\lambda_{j+n+1},\dots,\lambda_{m+n+p-2}}
\Big(a_{1},\dots,\\
&\quad a_{i-1}, f_{\lambda_{i},\dots,\lambda_{i+m-2}}
(a_{i},\dots, a_{i+m-2}), a_{i+m},\dots,a_{j-1},\\
&\quad a_{j}\oo{\lambda_{j}}g_{\lambda_{j+1},\dots,\lambda_{j+n-1}}(a_{j+1},
\dots, a_{j+n}), a_{j+n+1},\dots, a_{m+n+p-1}\Big),
\end{align*}
and
\begin{align*}
H_{i,j}''=&(-1)^{j+n-1}h_{\lambda_{1},\dots,\lambda_{i-1},
\lambda_{i}+\cdots+\lambda_{i+m-1},\lambda_{i+m},\dots,\lambda_{j-1},
\lambda_{j}+\cdots+\lambda_{j+n},\lambda_{j+n+1},\dots,\lambda_{m+n+p-2}}
\Big(a_{1},\dots,\\
&\quad a_{i-1},f_{\lambda_{i},\dots,\lambda_{i+m-2}}
(a_{i},\dots, a_{i+m-1}), a_{i+m},\dots, a_{j-1}, \\
&\quad g_{\lambda_{j},\dots,\lambda_{j+n-2}}(a_{j+1},
\dots, a_{j+n-1})\oo{\lambda_{j}+\cdots+\lambda_{i+n-1}}
a_{j+n}, a_{j+n+1},\dots, a_{m+n+p-1}\Big) \\
&+\sum_{q=j+n}^{m+n+p-2}(-1)^{q}h_{\lambda_{1},\dots,\lambda_{i-1},
\lambda_{i}+\cdots+\lambda_{i+m-1}, \lambda_{i+m},\dots,\lambda_{j-1},
\lambda_{j}+\cdots+\lambda_{j+n-1}, \lambda_{j+n},\dots,\lambda_{q-1},
\lambda_{q}+\lambda_{q+1},\lambda_{q+2},\dots,\lambda_{m+n+p-2}}\\
&\quad \Big(a_{1},\dots,a_{i-1}, f_{\lambda_{i},\dots,\lambda_{i+m-2}}
(a_{i},\dots, a_{i+m-1}), a_{i+m},\dots,a_{j-1},
\end{align*}
\begin{align*}
&\quad g_{\lambda_{j},\dots,\lambda_{j+n-2}}(a_{j},\dots, a_{j+n-1}),
a_{j+n},\dots, a_{q-1}, a_{q}\oo{\lambda_{q}}a_{q+1},a_{q+2},\dots,
a_{m+n+p-1}\Big)\\
&+(-1)^{m+n+p-1}h_{\lambda_{1},\dots,\lambda_{i-1},
\lambda_{i}+\cdots+\lambda_{i+m-1}, \lambda_{i+m},\dots,\lambda_{j-1},
\lambda_{j}+\cdots+\lambda_{j+n-1},\lambda_{j+n},\dots,\lambda_{m+n+p-3}}
\Big(a_{1},\dots,\\
&\quad a_{i-1}, f_{\lambda_{i},\dots,\lambda_{i+m-2}}
(a_{i},\dots, a_{i+m-1}), a_{i+m},\dots,a_{j-1},g_{\lambda_{j},\dots,
\lambda_{j+n-2}}(a_{j},\dots, a_{j+n-1}), \\
&\quad a_{j+n},\dots,a_{m+n+p-2}\Big)\oo{\lambda_{1}
+\cdots+\lambda_{m+n+p-2}}a_{m+n+p-1},
\end{align*}
for any $a_{1},\dots,a_{m+n+p-1}\in A$.
Then one can check that the following lemma is true.

\begin{lem}\label{lem:relation}
Let $A$ be an associative conformal algebra, and $f\in\mathcal{Z}^{m}(A, A)$,
$g\in\mathcal{Z}^{n}(A, A)$, $h\in\mathcal{C}^{p}(A, A)$, $m, n, p\geq1$.
Then for any $1\leq i\leq p-1$ and $m+i\leq j\leq m+p-1$, we have
\begin{equation*}
d_{m+n+p-1}((h\bullet_{i-1}f)\bullet_{j-1}g)_{\lambda_{1},\dots,\lambda_{m+n+p-2}}
(a_{1},\dots, a_{m+n+p-1})=H_{i,j}+H_{i,j}'+H_{i,j}''.
\end{equation*}
\end{lem}

For any $f\in\mathcal{Z}^{m}(A, A)$, $g\in\mathcal{Z}^{n}(A, A)$,
$h\in\mathcal{Z}^{1}(A, A)=\Der(A, A)$, $m, n\geq1$, we have
\begin{align*}
&\;(h\bullet(f\sqcup g))_{\lambda_{1},\dots,\lambda_{m+n-1}}(a_{1},\dots, a_{m+n})\\
=&\; h\Big(f_{\lambda_{1},\dots,\lambda_{m-1}}(a_{1},\dots, a_{m})
\oo{\lambda_{1}+\cdots+\lambda_{m}}g_{\lambda_{m+1},\dots,\lambda_{m+n-1}}
(a_{m+1},\dots, a_{m+n})\Big),\\
&\;((h\bullet f)\sqcup g)_{\lambda_{1},\dots,\lambda_{m+n-1}}(a_{1},\dots, a_{m+n})\\
=&\; h\Big(f_{\lambda_{1},\dots,\lambda_{m-1}}(a_{1},\dots, a_{m})\Big)
\oo{\lambda_{1}+\cdots+\lambda_{m}}g_{\lambda_{m+1},\dots,\lambda_{m+n-1}}
(a_{m+1},\dots, a_{m+n}),\\
&\;(f\sqcup (h\bullet g))_{\lambda_{1},\dots,\lambda_{m+n-1}}(a_{1},\dots, a_{m+n})\\
=&\; f_{\lambda_{1},\dots,\lambda_{m-1}}(a_{1},\dots, a_{m})
\oo{\lambda_{1}+\cdots+\lambda_{m}}h\Big(g_{\lambda_{m+1},\dots,\lambda_{m+n-1}}
(a_{m+1},\dots, a_{m+n})\Big).
\end{align*}
That is, $h\bullet(f\sqcup g)-(h\bullet f)\sqcup g-f\sqcup (h\bullet g)=0$
since $h\in\Der(A, A)$. But in general, for $h\in\mathcal{Z}^{p}(A, A)$,
$p\geq2$, the left hand side of the above equality may not be zero.
The following proposition show that this is given by a coboundary.

\begin{pro}\label{pro:coboundary}
Let $A$ be an associative conformal algebra, and $f\in\mathcal{Z}^{m}(A, A)$,
$g\in\mathcal{Z}^{n}(A, A)$, $h\in\mathcal{Z}^{p}(A, A)$, $m, n\geq1$ and $p\geq2$.
Set
\begin{equation*}
\mathcal{H}=\sum_{i=0}^{p-2}\sum_{j=m+i}^{m+p-2}(-1)^{(m-1)i+(n-1)j}
(h\bullet_{i}f)\bullet_{j}g.
\end{equation*}
Then $\mathcal{H}\in\mathcal{C}^{m+n+p-2}(A, A)$ and
\begin{equation*}
d_{m+n+p-2}(\mathcal{H})=(-1)^{(m-1)n}\Big(h\bullet(f\sqcup g)
-(-1)^{n(p-1)}(h\bullet f)\sqcup g-f\sqcup (h\bullet g)\Big).
\end{equation*}
\end{pro}

\begin{proof}
By Lemma \ref{lem:relation}, we have
\begin{equation*}
d_{m+n+p-1}((h\bullet_{i-1}f)\bullet_{j-1}g)_{\lambda_{1},\dots,\lambda_{m+n+p-2}}
(a_{1},\dots, a_{m+n+p-1})=H_{i,j}+H_{i,j}'+H_{i,j}''.
\end{equation*}
for any $1\leq i\leq p-1$, $m+i\leq j\leq m+p-1$ and $a_{1},\dots, a_{m+n+p-1}\in A$.
Note that if $m+i+1\leq j\leq m+p-2$,
\begin{align*}
&\;H_{i,j}+(-1)^{m-1}H_{i+1,j}'+(-1)^{m+n}H_{i+1,j+1}''\\
=&\;d_{p}(h)_{\lambda_{1},\dots,\lambda_{i}, \lambda_{i+1}+\cdots+\lambda_{i+m},
\lambda_{i+m+1},\dots,\lambda_{j}, \lambda_{j+1}+\cdots+\lambda_{j+n},
\lambda_{j+n+1},\dots,\lambda_{m+n+p-2}}\Big(a_{1},\dots, a_{i}, f_{\lambda_{i+1},
\dots,\lambda_{i+m-1}}(a_{i+1},\dots, \\
&\qquad a_{i+m}), a_{i+m+1},\dots, a_{j}, g_{\lambda_{j+1},\dots,\lambda_{j+n-1}}
(a_{j+1},\dots, a_{j+n}), a_{j+n+1},\dots, a_{m+n+p-1}\Big)\\
=&\;0,
\end{align*}
since $h\in\mathcal{Z}^{p}(A, A)$. In order for this equation to hold for any
$m+i\leq j\leq m+p-1$, we need the following signs.
\begin{itemize}
\item[] For $j=m, m+1,\dots, m+p+1$,
\begin{equation}\label{equ:re1}
H_{0,j}=(f\sqcup(h\bullet_{j-m}g))_{\lambda_{1},\dots,\lambda_{m+n+p-2}}
(a_{1},\dots, a_{m+n+p-1}).
\end{equation}
\item[] For $i=1, 2,\dots, p$,
\begin{equation}\label{equ:re2}
H'_{i,m+i-1}=(-1)^{m+i-1}(h\bullet_{j-1}(f\sqcup g))_{\lambda_{1},\dots,\lambda_{m+n+p-2}}
(a_{1},\dots, a_{m+n+p-1}).
\end{equation}
\item[] For $i=1, 2,\dots, p$,
\begin{equation}\label{equ:re3}
H''_{i,m+p}=(-1)^{m+n+p-1}((h\bullet_{i-1}f)\sqcup g)_{\lambda_{1},\dots,\lambda_{m+n+p-2}}
(a_{1},\dots, a_{m+n+p-1}).
\end{equation}
\end{itemize}
Then one can check that
\begin{equation}\label{equ:11}
\sum_{i=0}^{p-1}\sum_{j=m+i}^{m+p-1}(-1)^{(m-1)(i-1)+(n-1)(j-1)}
\Big(H_{i,j}+(-1)^{m-1}H'_{i+1,j}+(-1)^{m+n}H''_{i+1,j+1}\Big)=0.
\end{equation}
Moreover, by Lemma \ref{lem:relation}, one can check that
\begin{equation*}
d_{m+n+p-2}(\mathcal{H})_{\lambda_{1},\dots,\lambda_{m+n+p-2}}(a_{1},\dots, a_{m+n+p-1})
=\sum_{i=0}^{p-1}\sum_{j=m+i}^{m+p-1}(-1)^{(m-1)(i-1)+(n-1)(j-1)}
\Big(H_{i,j}+H'_{i+1,j}+H''_{i+1,j+1}\Big).
\end{equation*}
Thus, by equation (\ref{equ:11}), we obtain
\begin{align*}
0=&d_{m+n+p-2}(\mathcal{H})_{\lambda_{1},\dots,\lambda_{m+n+p-2}}
(a_{1},\dots, a_{m+n+p-1})+\sum_{i=0}^{p-1}
(-1)^{(m-1)(i-1)+(n-1)(m+p-2)+m+n}H''_{i+1,m+p}\\
&\qquad+\sum_{i=0}^{p-1}(-1)^{(m-1)(i-1)+(n-1)(m+i-1)+m-1}H'_{i+1,m+i}
+\sum_{j=m}^{m+p+1}(-1)^{(n-1)(j-1)-(m-1)}H_{0,j}.
\end{align*}
Bring equations (\ref{equ:re1}), (\ref{equ:re2}) and (\ref{equ:re3}) into the
above equation, we get
\begin{equation*}
d_{m+n+p-2}(\mathcal{H})=(-1)^{(m-1)n}f\sqcup (h\bullet g)
-(-1)^{(m-1)n}h\bullet(f\sqcup g)-(-1)^{(m-1)n+n(p-1)}(h\bullet f)\sqcup g.
\end{equation*}
Hence, we obtain the proposition.
\end{proof}

On the other hand, we have the following proposition.

\begin{pro}\label{pro:equation}
Let $A$ be an associative conformal algebra. For any $f\in\mathcal{C}^{m}(A, A)$,
$g\in\mathcal{C}^{n}(A, A)$, $h\in\mathcal{C}^{p}(A, A)$, $m, n, p\geq1$, we have
\begin{equation*}
(f\sqcup g)\bullet h=(f\bullet h)\sqcup g+(-1)^{m(p-1)}f\sqcup(g\bullet h).
\end{equation*}
\end{pro}

\begin{proof}
For $a_{1},\dots, a_{m+n+p-1}\in A$, we have
\begin{align*}
&\;((f\sqcup g)\bullet h)_{\lambda_{1},\dots,\lambda_{m+n+p-2}}(a_{1},\dots, a_{m+n+p-1}) \\
=&\;\sum_{i=1}^{m+n}(-1)^{(p-1)(i-1)}(f\sqcup g)_{\lambda_{1},\dots,\lambda_{i-1},
\lambda_{i}+\cdots+\lambda_{i+p-1},\lambda_{i+p},\dots,\lambda_{m+n+p-2}}
(a_{1},\dots, a_{i-1},\\
&\qquad h_{\lambda_{i},\dots,\lambda_{i+p-2}}(a_{i},\dots,
a_{i+p-1}), a_{i+p},\dots, a_{m+n+p-1}) \\
=&\;\sum_{i=1}^{m}(-1)^{(p-1)(i-1)}f_{\lambda_{1},\dots,\lambda_{i-1},
\lambda_{i}+\cdots+\lambda_{i+p-1},\lambda_{i+p},\dots,\lambda_{m+p-2}}
\Big(a_{1},\dots, a_{i-1},h_{\lambda_{i},\dots,\lambda_{i+p-2}}(a_{i},\dots
a_{i+p-1}),\\
&\qquad a_{i+p},\dots,a_{m+p-1}\Big)\oo{\lambda_{1}+\cdots+\lambda_{m+p-1}}
g_{\lambda_{m+p},\dots,\lambda_{m+n+p-2}}(a_{m+p},\dots, a_{m+n+p-1})\\
&\;+\sum_{i=m+1}^{m+n}(-1)^{(p-1)(i-1)}f_{\lambda_{1},\dots,\lambda_{m-1}}
(a_{1},\dots, a_{m})\oo{\lambda_{1}+\cdots+\lambda_{m}}g_{\lambda_{m+1},
\dots,\lambda_{i-1}, \lambda_{i}+\cdots+\lambda_{i+p-1}, \lambda_{i+p},\dots,
\lambda_{m+n+p-2}}\\
&\qquad \Big(a_{m+1},\dots, a_{i-1}, h_{\lambda_{i}+\cdots+\lambda_{i+p-2}}
(a_{i},\dots, a_{i+p-1}), a_{i+p}, \dots, a_{m+n+p-1} \Big) \\
=&\;(f\bullet h)_{\lambda_{1},\dots,\lambda_{m+p-2}}(a_{1}, \dots, a_{m+p-1})
\oo{\lambda_{1}+\cdots+\lambda_{m+p-1}}g_{\lambda_{m+p},\dots,\lambda_{m+n+p-2}}
(a_{m+p},\dots, a_{m+n+p-1})\\
&\;+(-1)^{m(p-1)}f_{\lambda_{1},\dots,\lambda_{m-1}}(a_{1},\dots, a_{m})
\oo{\lambda_{1}+\cdots+\lambda_{m}}(g\bullet h)_{\lambda_{m+1},
\dots,\lambda_{m+n+p-2}}(a_{m+1},\dots, a_{m+n+p-1}) \\
=&\;\Big((f\bullet h )\sqcup g+(-1)^{m(p-1)}f\sqcup(g\bullet h)
\Big)_{\lambda_{1},\dots,\lambda_{m+n+p-2}}(a_{1},\dots, a_{m+n+p-1}).
\end{align*}
Thus, $(f\sqcup g)\bullet h=(f\bullet h)\sqcup g+(-1)^{m(p-1)}f\sqcup(g\bullet h)$.
\end{proof}

Now, by Propositions \ref{pro:coboundary} and \ref{pro:equation}, we
obtain the Gerstenhaber algebra structure on the Hochschild cohomology
$\bigoplus_{i\geq1}\HH^{i}(A)$.

\begin{thm}\label{thm:Gers-alg}
Let $A$ be an associative conformal algebra. Then
$\Big(\bigoplus_{i\geq1}\HH^{i}(A),\; \sqcup,\; [-, -]\Big)$ is a
Gerstenhaber algebra.
\end{thm}

\begin{proof}
Let $f\in\mathcal{Z}^{m}(A, A)$, $g\in\mathcal{Z}^{n}(A, A)$,
$h\in\mathcal{Z}^{p}(A, A)$. If $m=n=p=0$, it is easy to see that
$[f\sqcup g, h]-[f, h]\sqcup g-(-1)^{m(p-1)}f\sqcup[g, h]=0$. If $n=p=0$ and $m\geq2$,
then $[g, h]=0$ and
\begin{align*}
&\;[f\sqcup g, h]_{\lambda_{1}, \dots, \lambda_{m-2}}(a_{1}, \dots, a_{m-1})\\
=&\;\sum_{i=1}^{m-1}(-1)^{i-1}f_{\lambda_{1}, \dots, \lambda_{i-1}, 0, \dots, \lambda_{i},
\dots,\lambda_{m-2}}(a_{1}, \dots, a_{i-1}, h, a_{i}, \dots, a_{m-1})\oo{-\partial}g\\
&\qquad\qquad\qquad\qquad +(-1)^{m-1}f_{\lambda_{1}, \dots, \lambda_{m-2}, -\partial}
(a_{1}, \dots, a_{m-1}, h)\oo{-\partial}g\\
=&\; ([f, h]\sqcup g)_{\lambda_{1}, \dots, \lambda_{m-2}}(a_{1}, \dots, a_{m-1}),
\end{align*}
for any $a_{1}, \dots, a_{m-1}\in A$. Thus, we also have the equation
$[f\sqcup g, h]-[f, h]\sqcup g-(-1)^{m(p-1)}f\sqcup[g, h]=0$.
Similarly, we get the equation whenever at least one of $m$, $n$, $p$ is zero, or $p=1$
by direct calculation.
If $p\geq2$, by Propositions \ref{pro:coboundary} and \ref{pro:equation}, we have
\begin{align*}
&\;[f\sqcup g, h]-[f, h]\sqcup g-(-1)^{m(p-1)}f\sqcup[g, h] \\
=&\;(f\sqcup g)\bullet h-(-1)^{(p-1)(m+n-1)}h\bullet(f\sqcup g)-(f\bullet h)\sqcup g
+(-1)^{(m-1)(p-1)}(h\bullet f)\sqcup g \\
&\;-(-1)^{m(p-1)}f\sqcup(g\bullet h)
+(-1)^{m(p-1)}(-1)^{(n-1)(p-1)}f\sqcup(h\bullet g) \\
=&\;-(-1)^{m(p-1)}(-1)^{(n-1)(p-1)}h\bullet(f\sqcup g)+(-1)^{p}(h\bullet f)\sqcup g
-(-1)^{(n-1)(p-1)}f\sqcup(h\bullet g)\\
=&\;\pm\Big(h\bullet(f\sqcup g)-f\sqcup(h\bullet g)
+(-1)^{np-n+1}(h\bullet f)\sqcup g \Big) \\
=&\;\pm d_{m+n+p-2}(\mathcal{H}).
\end{align*}
Thus, for any $f\in\HH^{m}(A)$, $g\in\HH^{n}(A)$, $h\in\HH^{p}(A)$, $m, n, p\geq0$,
\begin{equation*}
[f\sqcup g, h]-[f, h]\sqcup g-(-1)^{m(p-1)}f\sqcup[g, h]=0.
\end{equation*}
By the definition of Gerstenhaber algebra, we get this theorem.
\end{proof}


\section{Cohomology of split extension}\label{sec:exten}

In this section, we consider the relation between Hochschild cohomology of an
associative conformal algebra $A$ and Hochschild cohomology of it's split extension.
First, we give the definition of split extension of a associative conformal algebra.

\begin{defi}
Let $A$ be an associative conformal algebra, and $M$ be a conformal $A$-bimodule
with conformal product $\cdot\oo{\lambda}\cdot: M\times M\rightarrow M[\lambda]$. The
{\rm split extension} of $A$ by $M$ is an associative conformal algebra on $A\oplus M$,
where the conformal product is given by
\begin{equation*}
(a,\, u)\oo{\lambda}(b,\, v)=(a\oo{\lambda}b,\; a\ool{\lambda}v+u\oor{\lambda}b+
u\oo{\lambda}v),
\end{equation*}
for any $a, b\in A$ and $u, v\in M$. We denote by this associative conformal algebra
by $A\widehat{\oplus}M$.
\end{defi}

\begin{ex}
Let $A$ be an associative conformal algebra, and $M$ be a conformal $A$-bimodule.
The {\rm trivial extension} conformal algebra is an associative conformal algebra on
$A\oplus M$ with conformal product
\begin{equation*}
(a,\, u)\oo{\lambda}(b,\, v)=(a\oo{\lambda}b,\; a\ool{\lambda}v+u\oor{\lambda}b).
\end{equation*}
Trivial extension conformal algebras is a class of important split extension
conformal algebras.
\end{ex}

Given a split extension $A\widehat{\oplus}M$, there is an exact sequence of
$\cb[\partial]$-modules:
\begin{equation*}
\xymatrix@C=0.5cm{
0 \ar[r] & M\ar[rr]^{i\quad} && A\widehat{\oplus}M\ar[rr]<0.2ex>^{\quad p}
&& A\ar[r]\ar[ll]<0.5ex>^{\quad q} & 0 },
\end{equation*}
in which $p: (a, u)\mapsto a$, $i: u\mapsto (0, u)$, $q: a\mapsto (a, 0)$.
Clearly, $p$ and $q$ are morphisms of associative conformal algebras,
$p\circ q=\id_{A}$, $p$ induces a linear map $\tilde{p}: A\widehat{\oplus}M/
\partial(A\widehat{\oplus}M)\rightarrow A/\partial A$
and $q$ induces a linear map $\tilde{q}: A/\partial A\rightarrow
A\widehat{\oplus}M/\partial(A\widehat{\oplus}M)$.
Moreover, this sequence is a split sequence as conformal bimodule over $A$,
and a sequence of bimodule over $A\widehat{\oplus}M$, where the
$A\widehat{\oplus}M$-bimodule structure on $A$ given by the map $p$. Let $\varphi:
A\widehat{\oplus}M\rightarrow A\widehat{\oplus}M$ be a $\cb[\partial]$-module
homomorphism. Then $p\circ\varphi\circ q: A\rightarrow A$ is also a
$\cb[\partial]$-module homomorphism. For any $\varphi\in\mathcal{C}^{n}
(A\widehat{\oplus}M, A\widehat{\oplus}M)$, one can check that
$p\circ\varphi\circ q^{\otimes n}$ is also conformal sesquilinear, where
\begin{equation*}
q^{\otimes n}(a_{1},\dots,a_{n})=(q(a_{1}),\dots,q(a_{n})),
\end{equation*}
and $p$ is extended canonically to a $\cb[\partial]$-module homomorphism from
$A\widehat{\oplus}M[\lambda_{1},\dots,\lambda_{n-1}]$ to
$A[\lambda_{1},\dots, \lambda_{n-1}]$, i.e., $p(\sum x_{i_{1},\dots,i_{n-1}}
\lambda_{1}^{i_{1}}\cdots\lambda_{n-1}^{i_{n-1}})
=\sum p(x_{i_{1},\dots,i_{n-1}})\lambda_{1}^{i_{1}}\cdots\lambda_{n-1}^{i_{n-1}}$.

\begin{lem}\label{lem:map}
Let $A\widehat{\oplus}M$ be the split extension of $A$ by bimodule $M$. Let
$\tilde{d}_{n}$ (resp. $d_{n}$) denote the $n$-th differential in the Hochschild
complex of $A\widehat{\oplus}M$ with coefficients in $B$ (resp. of $A$ with
coefficients in $A$). Then we have
\begin{equation}\label{eq:0th}
d_{0}(\tilde{p}((a, u)+\partial(A\widehat{\oplus}M)))
=\tilde{p}\circ\tilde{d}_{0}((a, u)+\partial(A\widehat{\oplus}M))\circ\tilde{q}
\end{equation}
for any $(a, u)+\partial(A\widehat{\oplus}M)\in
A\widehat{\oplus}M/\partial(A\widehat{\oplus}M)$;
for $n\geq1$, and any $\varphi\in\mathcal{C}^{n}(A\widehat{\oplus}M,
A\widehat{\oplus}M)$,
\begin{equation}\label{eq:nth}
d_{n}(p\circ \varphi\circ q^{\otimes n})
=p\circ \tilde{d}_{n}(\varphi)\circ q^{\otimes (n+1)}.
\end{equation}
\end{lem}

\begin{proof}
First, if $n=0$, for any $b+\partial A\in A/\partial A$,
\begin{align*}
\tilde{p}\circ\tilde{d}_{0}((a, u)+\partial(A\widehat{\oplus}M))\circ\tilde{q}
(b+\partial A)&=\tilde{p}\circ\tilde{d}_{0}((a, u)+\partial(A\widehat{\oplus}M))
((b,0)+\partial(A\widehat{\oplus}M))\\
&=\tilde{p}((b,0)\oo{-\partial}(a, u)-(a, u)\oo{0}(b,0))\\
&=b\oo{-\partial}a-a\oo{0}b\\
&=d_{0}(\tilde{p}((a, u)+\partial(A\widehat{\oplus}M)))(b+\partial A).
\end{align*}
Second, for any $n\geq1$, $a_{1},\dots,a_{n+1}\in A$, since $p, q$ are
morphisms of associative conformal algebra and $p\circ q=\id_{A}$, we have
\begin{align*}
&\;d_{n}(p\circ \varphi\circ q^{\otimes n})_{\lambda_{1},\dots,\lambda_{n}}
(a_{1},\dots,a_{n+1})\\
=&\;a_{1}\oo{\lambda_{1}}(p\circ \varphi\circ q^{\otimes n})_{\lambda_{2},\dots,
\lambda_{n}}(a_{2},\dots,a_{n+1})\\
&\;+\sum_{i=1}^{n}(-1)^{i}(p\circ \varphi\circ q^{\otimes n})_{\lambda_{1},
\dots,\lambda_{i-1},\lambda_{i}+\lambda_{i+1}, \lambda_{i+2}\dots,\lambda_{n}}
(a_{1},\dots,a_{i-1}, a_{i}\oo{\lambda_{i}}a_{i+1}, a_{i+2},\dots,a_{n+1})\\
&\;+(-1)^{n+1}(p\circ \varphi\circ q^{\otimes n})_{\lambda_{1},\dots,
\lambda_{n-1}}(a_{1},\dots,a_{n})\oo{\lambda_{1}+\cdots+\lambda_{n}}a_{n+1}\\
=&\;p\circ q(a_{1})\oo{\lambda_{1}}p\circ \varphi_{\lambda_{2},\dots,
\lambda_{n}}(q(a_{2}),\dots,q(a_{n+1}))\\
&\;+\sum_{i=1}^{n}(-1)^{i}p\circ\varphi_{\lambda_{1},\dots,\lambda_{i-1},
\lambda_{i}+\lambda_{i+1},\lambda_{i+2}\dots,\lambda_{n}}(q(a_{1}),\dots,
q(a_{i-1}),q(a_{i}\oo{\lambda_{i}}a_{i+1}), q(a_{i+2}),\dots,q(a_{n+1}))\\
&\;+(-1)^{n+1}p\circ \varphi_{\lambda_{1},\dots, \lambda_{n-1}}(p(a_{1}),
\dots,p(a_{n}))\oo{\lambda_{1}+\cdots+\lambda_{n}}p\circ q(a_{n+1})\\
=&\;(p\circ \tilde{d}_{n}(\varphi)\circ q^{\otimes (n+1)})_{\lambda_{1},\dots,
\lambda_{n}}(a_{1},\dots,a_{n+1}).
\end{align*}
Thus $d_{n}(p\circ \varphi\circ q^{\otimes n})
=p\circ \tilde{d}_{n}(\varphi)\circ q^{\otimes (n+1)}$.
\end{proof}

For $n\geq0$, and any $\varphi\in\mathcal{C}^{n}(A, A)$, we denote its image
under the natural surjection $\mathcal{C}^{n}(A, A)\rightarrow \HH^{n}(A)$
by $[\varphi]$. Then by equation (\ref{eq:0th}), we have a linear map
$\Phi^{0}: \HH^{0}(A\widehat{\oplus}M)\rightarrow \HH^{0}(A)$,
$[(a, u)]\mapsto[\tilde{p}(a, u)]=[a]$. More general, we have the following corollary.

\begin{cor}\label{cor:HH-map}
Let $A\widehat{\oplus}M$ be the split extension of $A$ by bimodule $M$.
Then there exists a linear map $\Phi^{n}: \HH^{n}(A\widehat{\oplus}M)\rightarrow
\HH^{n}(A)$ given by $[\varphi]\mapsto [p\circ\varphi\circ q^{\otimes n}]$,
for each $n\geq1$.
\end{cor}

\begin{proof}
Assume that $\varphi\in\mathcal{C}^{n}(A\widehat{\oplus}M, A\widehat{\oplus}M)$
is a cocycle, that is, $\tilde{d}_{n}(\varphi)=0$. Then the formula
(\ref{eq:nth}) shows that $d_{n}(p\circ\varphi\circ q^{\otimes n})=0$,
and so that $p\circ\varphi\circ q^{\otimes n}$ is a cocycle.
If $\varphi\in\mathcal{C}^{n}(A\widehat{\oplus}M, A\widehat{\oplus}M)$ is a
coboundary, that is, $\varphi=\tilde{d}_{n-1}(\psi)$ for some $\psi\in
\mathcal{C}^{n-1}(A\widehat{\oplus}M, A\widehat{\oplus}M)$, then by formula
(\ref{eq:nth}), $p\circ\varphi\circ q^{\otimes n}=p\circ\tilde{d}_{n-1}(\psi)
\circ q^{\otimes n}=d_{n-1}(p\circ\psi\circ q^{\otimes (n-1)})$.
That is, $p\circ\varphi\circ q^{\otimes n}$ is also a coboundary.
Thus $\Phi^{n}$ is well-defined. Finally, it is easy to see that
$\Phi^{n}$ is $\cb$-linear.
\end{proof}

Now we can give the main conclusions of this section.

\begin{thm}\label{thm:alg-map}
Considering $\HH^{\ast}(A\widehat{\oplus}M)=\bigoplus_{i\geq 0}
\HH^{i}(A\widehat{\oplus}M)$ and $\HH^{\ast}(A)=\bigoplus_{i\geq 0}\HH^{i}(A)$
as rings under the cup product, the maps $\Phi^{n}$ induce a ring morphism
\begin{equation*}
\Phi^{\ast}: \HH^{\ast}(A\widehat{\oplus}M)\rightarrow\HH^{\ast}(A).
\end{equation*}
\end{thm}

\begin{proof}
Let $\eta=[f]\in\HH^{s}(A\widehat{\oplus}M)$ and
$\zeta=[g]\in\HH^{t}(A\widehat{\oplus}M)$,
where $f\in\mathcal{C}^{s}(A\widehat{\oplus}M, A\widehat{\oplus}M)$
and $g\in\mathcal{C}^{t}(A\widehat{\oplus}M, A\widehat{\oplus}M)$.
First, if $s=t=0$, then $f=(a, u)$ and $g=(a', u')$ for some $a, a'\in A$ and $u, u'\in M$,
and
$$
\Phi^{0}(\eta)\sqcup\Phi^{0}(\zeta)=[a\oo{\lambda}a']
=[\tilde{p}((a, u)\oo{\lambda}(a', u'))]=\Phi^{0}(\eta\sqcup\zeta).
$$
Second, if $s, t\geq1$, then by Corollary \ref{cor:HH-map}, $\Phi^{s}(\eta)=[p\circ f\circ
q^{\otimes s}]$, $\Phi^{t}(\zeta)=[p\circ g\circ q^{\otimes t}]$, and
\begin{align*}
&\;\Big((p\circ f\circ q^{\otimes s})\sqcup(p\circ g\circ q^{\otimes t})
\Big)_{\lambda_{1},\dots,\lambda_{s+t-1}}(a_{1},\dots,a_{s+t})\\
=&\;(p\circ f\circ q^{\otimes s})_{\lambda_{1},\dots,\lambda_{s-1}}(a_{1},
\dots,a_{s})\oo{\lambda_{1}+\cdots+\lambda_{s}}(p\circ g\circ
q^{\otimes t})_{\lambda_{s+1},\dots,\lambda_{s+t-1}}(a_{s+1},\dots,a_{s+t})\\
=&\;p\circ f_{\lambda_{1},\dots,\lambda_{s-1}}(q(a_{1}),\dots,q(a_{s}))
\oo{\lambda_{1}+\cdots+\lambda_{s}}(p\circ g_{\lambda_{s+1},\dots,\lambda_{s+t-1}}
(q(a_{s+1}),\dots,q(a_{s+t}))\\
=&\;(p\circ (f\sqcup g))_{\lambda_{1},\dots,\lambda_{s+t-1}}(q(a_{1}),\dots,q(a_{s+t}))\\
=&\;(p\circ (f\sqcup g)\circ q^{\otimes (s+t)})_{\lambda_{1},\dots,\lambda_{s+t-1}}
(a_{1},\dots,a_{s+t}).
\end{align*}
Thus $(p\circ f\circ q^{\otimes s})\sqcup(p\circ g\circ q^{\otimes t})
=p\circ (f\sqcup g)\circ q^{\otimes (s+t)}$. Taking the cohomology classes,
we get $\Phi^{s}(\eta)\sqcup\Phi^{t}(\zeta)=\Phi^{s+t}(\eta\sqcup\zeta)$.
Similarly, if $s=0$ and $t\geq1$, or $s\geq1$ and $t=0$, we also have
$\Phi^{s}(\eta)\sqcup\Phi^{t}(\zeta)=\Phi^{s+t}(\eta\sqcup\zeta)$.
Therefore $\Phi^{\ast}: \HH^{\ast}(A\widehat{\oplus}M)\rightarrow\HH^{\ast}(A)$
is a ring morphism.
\end{proof}

The algebra homomorphism $\Phi^{\ast}$ is often called the Hochschild projective
morphism. For the Hochschild cohomology of associative algebras and the Poisson
cohomology of Poisson algebras, the Hochschild projective morphism
is not a morphism of graded Lie algebras in general. And in \cite{AGST} and
\cite{ZW}, the authors have given some conditions for the Hochschild projective
morphism to be surjective. For the associative conformal algebra,
we need to further consider the condition that $\Phi^{\ast}$ is surjective.
The sequence $\xymatrix@C=0.5cm{ 0 \ar[r] & M\ar[r] & A\widehat{\oplus}M
\ar[r]<0.2ex>& A\ar[r]\ar[l]<0.5ex> & 0 }$ consider in this section is an
abelian extension of $A$ by $M$. In the following work \cite{HZ}, we will give the
classification of the non-abelian extensions of associative conformal algebras by using
the Maurer-Cartan elements of a differential graded Lie algebra related to cohomology
of associative conformal algebras.

\bigskip
\noindent
{\bf Acknowledgements }
This work was financially supported by National
Natural Science Foundation of China (No. 11771122, 11801141, 12201182), China
Postdoctoral Science Foundation (2020M682272) and NSF of Henan Province
(212300410120). The authors are indebted to the referee for his/her help comments and
suggestions which have improved the article.


\begin{thebibliography}{abc}

\bibitem{AGST} I. Assem, M.A. Gatica, R. Schiffler, R. Taillefer,
             Hochschild cohomology of relation extension algebras,
             J. Pure Appl. Algebra {\bf 220}, 2016, 2471--2499.

\bibitem{BDK} B. Bakalov, A. D'Andrea, V.G. Kac,
              Theory of finite pseudoalgebras.
              Adv. Math. {\bf 162}, 2001, 1--140.

\bibitem{BFK} L.A. Bokut, Y. Fong, W.F. Ke,
              Free associative conformal algebras,
              Proc. of the 2nd Tainan-Moscow Alg. and Comb. Workshop, Tainan 1997,
              pp. 13--25.

\bibitem{BFK1}L.A. Bokut, Y. Fong, W.F. Ke,
              Grobner-Shirshov bases and composition lemma for associative conformal
              algebras: an example,
              Contemp. Math. {\bf 264}, 2000, 63--90.

\bibitem{BFK2}L.A. Bokut, Y. Fong, W.F. Ke,
              Composition-Diamond lemma for associative conformal algebras,
              J. Algebra {\bf 272}, 2004, 739--774.

\bibitem{BKV} B. Bakalov, V.G. Kac, A. Voronov,
              Cohomology of conformal algebras.
              Comm. Math. Phys. {\bf 200}, 1999, 561--589.

\bibitem{BPZ} A.A. Belavin, A.M. Polyakov, A.B. Zamolodchikov,
              Infinite conformal symmetry in two-dimensional quantum field theory,
              Nuclear Phys. {\bf 241}, 1984, 333--380.

\bibitem{CK}  S. Cheng and V.G. Kac,
              Conformal modules,
              Asian J. Math.   {\bf 1}, 1997, 181--193.

\bibitem{CKW} S. Cheng, V.G. Kac, M. Wakimoto,
              Extensions of conformal modules,
              Topological field theory, primitive forms and related topics (Kyoto),
              Progress in Math. 160. Birkhauser, Boston 1998, pp: 33--57.

\bibitem{DK}  A. D'Andrea, V.G. Kac,
              Structure theory of finite conformal algebras,
              Sel. Math., New Ser. {\bf 4}, 1998, 377--418.

\bibitem{Das}  A. Das,
              Gerstenhaber algebra structure on the cohomology of a
              Hom-associative algebra,
              Proc. Indian Acad. Sci.(Math. Sci.) {\bf 130}, 2020, Paper
              No.20, 20 pp.

\bibitem{Do}  I.A. Dolguntseva,
              The Hochschild cohomology for associative conformal algebras,
              Algebra Logic {\bf 46}, 2007, 373--384.

\bibitem{Do1} I.A. Dolguntseva,
              Triviality of the second cohomology group of the conformal algebras
              ${\rm Cend}_n$ and ${\rm Cur}_n$,
              St. Petersburg Math. J. {\bf 21}, 2010, 53--63.

\bibitem{Ger} M. Gerstenhaber,
             The cohomology structure of an associative ring,
             Ann. of Math. (2) 78 (1963), 267-288.

\bibitem{HZ} B. Hou, J. Zhao,
             On non-abelian extensions associative conformal algebras,
             in preparation.

\bibitem{Ka} V.G. Kac,
              Vertex algebras for beginners, 2nd Edition,
              Amer. Math. Soc. Providence, RI, 1998.

\bibitem{Ka1} V.G. Kac,
              Formal distribution algebras and conformal algebras,
              XII-th Int. Congr. Math. Phys.  Int. Press, Cambridge, MA 1999,
              pp. 80--97.

\bibitem{Ko}  P.S. Kolesnikov,
              Simple associative conformal algebras of linear growth,
              J. Algebra {\bf 295}, 2006, 247--268.

\bibitem{Ko1} P.S. Kolesnikov,
              Associative conformal algebras with finite faithful representation,
              Adv. Math. {\bf 202}, 2006, 602--637.

\bibitem{Ko2} P.S. Kolesnikov,
              On the Wedderburn principal theorem in conformal algebras,
              J. Algebra Appl. {\bf 6}, 2007, 119--134.

\bibitem{Ko3} P.S. Kolesnikov,
              On finite representations of conformal algebras.
              J. Algebra {\bf 331}, 2011, 169--193.

\bibitem{KK}  P.S. Kolesnikov, R.A. Kozlov,
              On the Hochschild cohomologies of associative conformal algebras with
              a finite faithful representation,
              Comm. Math. Phys. {\bf 369}, 2019, 351--370.

\bibitem{Koz} R.A. Kozlov,
              Hochshild cohomology of the associative conformal algebra
              ${\rm Cend}_{1,x}$,
              Algebra and Logic {\bf 58}, 2019, 36--47.

\bibitem{Lib} J.I. Liberati,
              Cohomology of associative $H$-pseudoalgebras,
              arXiv: 2012.13832.

\bibitem{Re}  A. Retakh,
              Associative conformal algebras of linear growth,
              J. Algebra {\bf 237}, 2001, 769--788.

\bibitem{Re1} A. Retakh,
              On associative conformal algebras of linear growth II,
              J. Algebra {\bf 304}, 2006, 543--556.

\bibitem{Ro}  M. Roitman,
              Universal enveloping conformal algebras,
              Sel. Math. New Ser. {\bf 6}, 2000, 319--345.

\bibitem{Ro1} M. Roitman,
              On embedding of Lie conformal algebras into associative
              conformal algebras,
              J. Lie theory {\bf 15}, 2005, 575-588.

\bibitem{ZW} C. Zhu, G. Wu,
             Poisson cohomology of trivial extension algebras,
             Bull. Iran. Math. Soc. {\bf 47}, 2021, 535-552.



\end{thebibliography}
 \end{document}